\documentstyle[a4,leqno,amsfonts,12pt]{article}
\newcommand{\C}{{\bf C}}

\newcommand{\OC}{\overline{{\bf C}}}
\newcommand{\R}{{\bf R}}
\newcommand{\He}{{\bf H}}
\newcommand{\N}{{\bf N}}

\newcommand{\bP}{{\bf P}}
\newcommand{\D}{{\bf D}}
\newcommand{\de}{\delta}

\newcommand{\mb}{\mbox}

\newcommand{\beq}{\begin{equation}}
\newcommand{\eeq}{\end{equation}}
\newcommand{\oge}{\succeq}
\newcommand{\ole}{\preceq}

\newcommand{\ov}{\overline}
\newcommand{\al}{\alpha}
\newcommand{\be}{\beta}
\newcommand{\Om}{\Omega}

\newcommand{\z}{\zeta}

\newcommand{\kap}{\mb{ cap}}
\newcommand{\ga}{\gamma}

\newcommand{\la}{\lambda}
\newcommand{\Ga}{\Gamma}

\newtheorem{th}{Theorem}

\newtheorem{lem}{Lemma}

\newcommand{\ueberschrift}{\bigskip\goodbreak\noindent\bigskip}
\newcounter{theabsatz}
\newcommand{\absatz}[1]{\stepcounter{theabsatz} \ueberschrift
                           {\large \bf \arabic{theabsatz}. {#1}} \setcounter{equation}{0}}

\begin{document}
\mathsurround=2pt

\begin{center}
{\bf \large ON CHEBYSHEV POLYNOMIALS in the COMPLEX PLANE}\\[2ex]

V. V. ANDRIEVSKII \\[2ex]

 Department of Mathematical
Sciences,
 Kent State University,\\
 Kent, OH 44242, U.S.A.\\
 email: andriyev@math.kent.edu\\[2ex]

\end{center}

 {\bf Abstract.}
The estimates of the uniform norm of the Chebyshev polynomials
 associated with a compact set $K$ in the complex plane are established.
 These estimates are exact (up to a constant factor) in the case
 where $K$ consists of  a finite number of quasiconformal curves or  arcs. 
The case where $K$
 is a uniformly perfect  subset of the real line
 is also studied.


\absatz{Introduction and main results}

Let $K\subset\C$ be a compact set in the complex plane $\C$
with a connected complement $\Om:=\ov{\C}\setminus K$, where
$\ov{\C}:=\C\cup\{\infty\}.$ We assume that $\kap(K)>0$,
where$\kap(K)$ denotes the logarithmic capacity of $K$
 (see \cite{pom} - \cite{saftot}).
Denote by $T_n(z)=T_n(z,K),n\in\N:=\{ 1,2,\ldots\}$
the $n$-th Chebyshev polynomial associated with $K$, i.e.,
$T_n(z)=z^n+c_{n-1}z^{n-1}+\ldots+c_0, c_k\in\C,$ is the (unique)
monic polynomial which minimizes the supremum norm $||T_n||_K:=
\sup_{z\in K}|T_n(z)|$ among all monic polynomials of the same
degree.

It is well-known (see, for example, \cite[Theorem 5.5.4 and
Corollary 5.5.5]{ran}) that
$$
||T_n||_K\ge\kap(K)^n\quad\mb{and}\quad
\lim_{n\to\infty}||T_n||_K^{1/n}=\kap(K).
$$
We  are interested in estimates  from above for the quantity
$$
t_n(K):=\frac{||T_n||_K}{\kap(K)^n}\, .
$$
We refer the reader to \cite{smileb}-\cite{sue}, \cite{tot12}-\cite{totvar}, \cite{wid}, 
\cite{csz}, \cite{and151}
and many references therein for a comprehensive survey of this subject.

First, let $K$ consist
 of disjoint  closed connected sets (continua) $K^j,j=1,2,\ldots,m$,
i.e.,
\beq\label{1.1}
K=\cup_{j=1}^mK^j;\quad K^j\cap K^k=\emptyset\mb{ for }j\neq k;
\quad\mb{ diam}(K^j)>0.
\eeq
Here
$$
\mb{diam}(S):=\sup_{z,\z\in S}|z-\z|,\quad S\subset\C.
$$
\begin{th}
\label{th1}
Under the above assumptions,
\beq\label{1.2}
t_n(K)= O(\log n)\quad \mb{for }n\to\infty.
\eeq
\end{th}
If more information is known about the geometry of $K$, (\ref{1.2})
can be improved, for example,  in the following
way.
A Jordan curve $L\subset\C$ is called quasiconformal
(see \cite{ahl} or \cite[p. 100]{lehvir}) if
for every $z_1,z_2\in L$,
\beq\label{1.3}
\mb{diam}(L(z_1,z_2))\le \beta_L|z_2-z_1|,
\eeq
where $L(z_1,z_2)$ is the smaller subarc of $L$ between $z_1$ and $z_2$,
 a constant $\beta_L>1$ depends only on $L$. Any subarc of a quasiconformal
curve is called a quasiconformal arc.
\begin{th}\label{th2} 
Let each $K^j$ in (\ref{1.1})
be either  a quasiconformal arc
or a closed Jordan domain bounded by  a quasiconformal curve. Then
\beq\label{1.4}
t_n(K)=O(1)\quad \mb{for }n\to\infty.
\eeq
\end{th}
The estimate (\ref{1.4})
was  proved by other methods in \cite{wid} and recently in \cite{tot14}
(for  sufficiently smooth $\partial K^j$), in \cite{totvar} (for piecewise
sufficiently smooth $\partial K^j$), and in \cite{and151} (for quasismooth 
in the sense of Lavrentiev
$\partial K^j$).

The question whether (\ref{1.4}) does hold for a general continuum 
seems to be still open. In the Oberwolfach meeting
(see \cite{jah} or \cite[p. 365]{pom72})
Pommerenke asked about an example for a continuum $K$ such that the
sequence
$\{t_n(K)\}$ is unbounded. It is mentioned in \cite[p. 365]{pom72} that
``D. Wrase in Karlsruhe has shown that an example constructed by J. Clunie
\cite{clu} for a different purpose has the required property".
But we could not find the proof of this result.

Moreover, in the case where
 $K$ is a continuum, one of the major sources
 for estimates of  $t_n(K)$ are Faber polynomials $F_n(z)=F_n(z,K)$
 associated with $K$ (see \cite{smileb}, \cite{sue}). Gaier
 \cite[Theorem 2]{gai99}, using the same example  by  Clunie,
\cite{clu} has shown that
 there exist a continuum $K_*$ bounded by a quasiconformal
 curve with $\kap(K_*)=1$, a positive constant $\al$ and an infinite set
 $\Lambda\subset\N$ such that for the (monic) polynomial $F_n(z)=F_n(z,K_*)$
  we have
 $$
 ||F_n||_{K_*}> n^\al,\quad n\in\Lambda .
 $$
 Note that the first result of this kind (without the restriction on $K_*$
 to be a quasidisk) was proved by Pommerenke \cite{pom66}.

Hence, Theorem \ref{th1} and
  Theorem \ref{th2} reveal the  
	essential difference between
	 the Chebyshev
  and the Faber polynomials. It is worth pointing out that
	the case of multiply connected $\Om$ presents a more delicate problem
	(see for example \cite{wid}).

Let now  $K\subset\R$, where $\R$ is the real line,
consist of an infinite number of components.
According to \cite[Theorem 4.4]{gon} in this case
$\{t_n(K)\}$ can increase  faster than any sequence $\{t_n\}$
satisfying
$t_n\ge1$ and  $\lim_{n\to\infty}(\log t_n)/n=0$. Therefore, in order to have particular
bounds for $t_n(K)$ some additional assumptions on $K$ are needed. We assume that
$K$ is uniformly perfect,
which according to Beardon and Pommerenke \cite{beapom} means that
there exists a  constant $0<\ga_K<1$ such that for  $z\in K$,
$$
 K\cap\{ \z:\, \ga_Kr\le|z-\z|\le r\}\neq\emptyset,\quad
0<r<\mb{diam}(K).
$$
The classical Cantor set is an example of a uniformly perfect set.
Pommerenke \cite{pom79} has shown that uniformly perfect sets can
be described using a density condition  in terms of the
logarithmic capacity. Namely, $K$ is uniformly perfect if and only if there
exists a constant $0<\la_K<1$ such that for  $z\in K$,
 \beq\label{1.5}
 \kap(K\cap\{\z:\, |\z-z|\le r\})\ge \la_K\, r,\quad 0<r<\mb{diam}(K).
 \eeq
Note that sets
satisfying (\ref{1.5}) play a significant role in the solution of
the inverse problem of the constructive theory of functions of a
complex variable. We refer to \cite{tam} where they are called
$c$-dense sets. Other interesting properties of the uniformly perfect sets can be found in
\cite[pp. 343--345]{garmar}.
\begin{th} \label{th3} For a uniformly perfect set $K\subset\R$  
there exists a constant $c=c(K)>0$ such that
\beq\label{1.6}
t_n(K)=O(n^c)\quad \mb{for }n\to\infty.
\eeq
\end{th}
Following Carleson \cite{car}
we say that a compact set $K\subset\R$ is homogeneous if there is a constant
$\eta_K>0$ such that for all $x\in K$,
$$
|K\cap(x-r,x+r)|\ge\eta_K r,\quad 0<r<\mb{diam}(K).
$$
Here, $|S|$ is the linear measure (length) of a
(Borel) set $S\subset\C$ (see \cite[p. 129]{pom}).
The Cantor sets of positive length are examples
of homogeneous sets  (see
\cite[p. 125]{pehyud}).
 Recently Christiansen, Simon, and Zinchenko \cite{csz} have 
shown that for the homogeneous subsets of the real line
the term $O(n^c)$ in (\ref{1.6}) can be replaced by $O(1)$. 
It is worth pointing out that there is a principal difference between the  above 
mentioned classes of compact sets, i.e., $K$ is the Parreau-Widom set in the case
of the homogeneous $K\subset \R$ and it is not, in general,
the Parreau-Widom set in the case
of the uniformly perfect $K$. See \cite{csz} for more details.

In what follows, 
we  use
the convention that $c,c_1,\ldots$
denote positive constants 
(different in different sections)
that are either absolute or they depend only
on $K$; otherwise, the dependence on  other parameters is explicitly
stated.
 For the nonnegative functions $a$
and $b$ we write $a\preceq b$ if $a\le cb$, and $a\asymp b$ if
$a\preceq b$ and $b\preceq a$ simultaneously.

We also use the additional notation
$$
d(z,S):=\mb{dist}(\{z\},S):=\inf_{\z\in S}|z-\z|,\quad
z\in\C, S\subset\C.
$$

\absatz{The basic potential-theoretic functions}

Let $K$ be as in (\ref{1.1}).
Following Widom \cite{wid}, we extend the concept of Faber polynomials to
the case of compact sets with the finite number of connected components.
Since in \cite{wid} all $\partial K^j$ are sufficiently smooth curves, we need
to add some purely technical details.
Denote by
$ g_{\Om}(z,z_0), \: z,z_0 \in \Omega,$
  the Green function for $\Omega$ with
 pole at $z_0$. It has a multiple-valued harmonic conjugate $\tilde{g}_{\Om}(z,z_0)$.
 Thus, the analytic function
    $$
    \Phi_{\Om}(z,z_0) :=  \exp(g_{\Om}(z,z_0) + i \tilde{g}_{\Om}(z,z_0))
    $$
is also multiple-valued.
We  write $g_{\Om}(z)$,
$\tilde{g}_{\Om}(z)$,  and $\Phi_{\Om}(z)$ in the case $z_0 = \infty$.

Let 
$$
g_\Om(z):=0,\quad z\in K,
$$
$$
K_s:=\{z\in\C:g(z)\le s\},\Om_s:=\OC\setminus K_s, \quad s>0.
$$
Then for $z\in\Om_s$,
$$
g_{\Om_s}(z)=g_{\Om}(z)-s, \Phi_{\Om_s}(z)=e^{-s}\Phi_{\Om}(z).
$$
For $n \in \N$, 
if $\Phi_{\Om}(z)^n$ is  single-valued
in $\Omega_{1/n^2}$, we set
   $$
   W_n(z) := \Phi_{\Om}(z)^n, \quad z \in \Omega_{1/n^2}.
  $$
If $\Phi_{\Om}(z)^n$ is multiple-valued in $\Omega_{1/n^2}$, then according to
\cite[pp. 159, 211]{wid}
 there exist $q \le m-1$ points $z_{1,n},\ldots,z_{q,n}\in\Om_{1/n^2}$  such that  
the function
   $$
   W_n (z) := \Phi_{\Om}(z)^n \prod^q_{l=1} \Phi_{\Om_{1/n^2}}(z,z_{l,n})^{-1}, \quad z \in 
	\Omega_{1/n^2},
   $$
 is single-valued in $\Omega_{1/n^2}$.
Moreover, all $z_{l,n}$ lie in the convex hull of $K_{1/n^2}$.

In both cases we consider the entire function
\beq\label{2n.1}
    F_n(z) := \frac{1}{2\pi i}  \int_{C_n}  \frac{W_n(\zeta)}
    {\zeta-z} \, d\zeta, \quad z\in\C,
   \eeq
 where $C_n\subset\Om_{1/n^2}\setminus\{\infty\}$ is a Jordan curve, 
oriented in the positive direction, containing
$K_{1/n^2}$ and $z$ in its interior.

 Since all points $z_{l,n}$ are in the convex hull of $K_1$, by 
the symmetry property of the Green function, we obtain
\begin{eqnarray*}
\left|\prod^q_{l=1}\Phi_{\Om_{1/n^2}}(\infty,z_{l,n})\right|
&=&\exp\left(\sum^q_{l=1}g_{\Om_{1/n^2}}(z_{l,n})\right)
\\
&\le&
\exp\left(\sum^q_{l=1}g_{\Om}(z_{l,n})\right)\ole 1.
\end{eqnarray*}
For $z\in\Om_{1/n^2}\setminus\{\infty\}$, let
$C_n'\subset\Om_{1/n^2}\setminus\{\infty\}$ be any Jordan curve, 
oriented in the positive direction, containing
$K_{1/n^2}$ in its interior and $z$ in its exterior.
Since by the Cauchy formula 
    $$
    F_n(z) = W_n(z) + \frac{1}{2\pi i} \; \int_{C_n'} \; \frac{W_n(\zeta)}
    {\zeta-z} d\zeta,
    $$
we see that $F_n(z)=\al_nz^n+\ldots$ is a polynomial with the property
 \begin{eqnarray}
 |\al_n|&=&
 \lim_{z\to\infty}\left|\frac{F_n(z)}{z^n}\right|=
\lim_{z\to\infty}\left|\frac{W_n(z)}{z^n}\right|
\nonumber\\
\label{2.2}
  &=&\kap(K)^{-n}\prod_{l=1}^q\left|\Phi_{\Om_{1/n^2}}(z_{l,n})\right|^{-1}
	\oge \kap(K)^{-n}.
 \end{eqnarray}
Now let $K$ consist of one component, i.e., $m=1$ and let $\Phi_\Om:\Om\to\D^*:=\{
w:|w|>1\}$ be the Riemann conformal mapping with
$\Phi_\Om(\infty)=\infty, \Phi_\Om'(\infty)>0$.
We follow a technique of \cite[Chapter IX]{dzj},
\cite[p. 387]{andbla} and for $k,n\in\N,k\ge2$, consider the Dzjadyk polynomial kernel
$$
K_{1,1,k,n}(\z,z),\quad \z\in \Om\setminus\{\infty\}
, z\in K,
$$
which is a polynomial with respect to $z$ of degree $(k+3)k(n-1)-1$
with continuous coefficients depending on $\z$.

According to \cite[p. 389, Theorem 2.4]{andbla} we have
\beq\label{2nn.1}
\left|\frac{1}{\z-z}-K_{1,1,k,n}(\z,z)\right|
\le c_1\frac{|\tilde{\z}_{1/n}-\z|^k}{|\z-z||\tilde{\z}_{1/n}-z|^k}
\left(1+\left|\frac{\z-z}{\tilde{\z}_{1/n}-z}\right|\right)^k,
\eeq
where $c_1=c_1(K,k)$ and
$$
\tilde{\z}_\de:=\Phi^{-1}_\Om((1+\de)\Phi_\Om(\z)),\quad \de>0.
$$
Let $w:=\Phi_\Om(\z),\tilde{w}_{1/n}:=w(1+1/n)$.
A straightforward calculation shows that
for $n>32^2$ and $|w|\ge1+32/n$, we have
$$
\frac{|\tilde{w}_{1/n}-w|}{|\tilde{w}_{1/n}|-1}<\frac{1}{32}.
$$
Therefore, \cite[p. 23, Lemma 2.3]{andbla} implies 
$$
\frac{|\tilde{\z}_{1/n}-\z|}{|\tilde{\z}_{1/n}-z|}\le 16
\frac{|\tilde{w}_{1/n}-w|}{|\tilde{w}_{1/n}|-1}<\frac{1}{2},
$$
i.e.,
$$
\left|\frac{1}{\z-z}-K_{1,1,k,n}(\z,z)\right|\le c_2
\frac{d(\z,K)^k}{|\z-z|^{k+1}},\quad c_2=c_2(K,k).
$$
We summarize our reasoning as follows. Given $k\in\N$, 
there exist sufficiently large constants $n_0=n_0(k)$ and
$c_3=c_3(k)$ such that
for any
integer $n>n_0$ and $\z$ with $|\Phi_\Om(\z)|-1\ge c_3/n$,
there exists a polynomial
$$
p_{n,k,\z,K}(z)=\sum_{l=0}^{n}a_{l,k,K}(\z)z^l,
$$
where $a_{l,k,K}$ are continuous functions of $\z$, satisfying
\beq\label{2.4}
\left|\frac{1}{\z-z}-p_{n,k,\z,K}(z)\right|\le c_2
\frac{d(\z,K)^k}{|\z-z|^{k+1}},\quad  z\in K.
\eeq
Indeed, to get (\ref{2.4}) we can take
$$
p_{n,k,\z,K}(z):=K_{1,1,k,N}(\z,z),\quad
N:=\left\lfloor \frac{n}{k(k+3)}\right\rfloor .
$$
Furthermore, by virtue of (\ref{2nn.1}), for $\z\in\Om$ with
$c_4\le|\Phi_\Om(\z)|-1\le c_5$ we have
\beq\label{2nn.2}
\left|\frac{1}{\z-z}-p_{n,k,\z,K}(z)\right|\le
\frac{c_6}{n^{k}},\quad  z\in K, c_6=c_6(c_4,c_5,K,k).
\eeq
Let $K$ now be as in (\ref{1.1}) with $m>1$. Denote by $r_K>0$ any fixed number  
such that $K_{r_K}$ consists of exactly $m$
components, i.e.,
$$
K_{r_K}=\cup_{j=1}^mK_{r_K}^j,\quad K^j\subset K_{r_K}^j.
$$
Let $\Om^j:=\OC\setminus K^j$.
The maximum principle for  the appropriate linear combination of
harmonic functions 
$g_\Om$ and $\log|\Phi_{\Om^j}|$ in $K_{r_K}^j\setminus K^j$ shows
that
\beq\label{2.6}
g_\Om(\z)\ole \log|\Phi_{\Om^j}(\z)|,\quad \z\in K_{r_K}^j\setminus K^j.
\eeq

For sufficiently large $v\in\N$, $\z\in K_{r_K}^j\setminus K^j$ 
with $|\Phi_{\Om^j}(\z)|-1\ge c_3/v$, and $z\in K^l,l=1,\ldots m$,
by virtue of (\ref{2.4}) and (\ref{2nn.2}), 
applied for the continuum $K^l$,
we have
\beq\label{2.7}
\left|\frac{1}{\z-z}-p_{v,k,\z,K^l}(z)\right|\le c_7
\left\{\begin{array}{ll}
   \displaystyle
\frac{d(\z,K^j)^k}{|\z-z|^{k+1}}
 &\mb{ if } l=j,
\\[2ex]
n^{-k}&\mb{ if }l\neq j.
\end{array}\right.
\eeq
Here $c_7:= c_2+c_6$.

For $\z$ as in (\ref{2.7}) and $l=1,\ldots, m,$
consider  functions
$$
h_l(z):=
\left\{\begin{array}{ll}
   \displaystyle
 1
 &\mb{ if } z\in K^l,
\\[2ex]
0&\mb{ if }z\in K\setminus K^l,
\end{array}\right.
$$
$$
f_{\z,l}(z):=\frac{h_l(z)}{\z-z},\quad z\in K,
$$
so that 
$$
\frac{1}{\z-z}=\sum_{l=1}^m f_{\z,l}(z),\quad z\in K.
$$
Since $h_l$ can be extended analytically to $K_{r_K}$,
by the Walsh approximation theorem \cite[pp. 75-76]{wal} there is
$u_0=u_0(K)\in\N$, such that for any integer $u>u_0$, there exists a polynomial
$q_{u,l}\in\bP_u$ satisfying
\beq\label{2.8}
||h_l-q_{u,l}||_K\le e^{-ur_K}.
\eeq
For sufficiently large $n$ and 
$\z\in K_{r_K}^j\setminus K^j$ 
with $|\Phi_{\Om^j}(\z)|-1\ge c^*/n\ge c_3/v$,
where  $v$ and the constant $c^*$ are to be chosen later,
consider the polynomial 
$$
t_{u,v,k,\z,l}:=q_{u,l}p_{v,k,\z,K^l}\in\bP_{u+v}.
$$
Let
$$
R_K:=\max_{1\le j\le m}||\log|\Phi_{\Om^j}|||_{\partial K}.
$$
Since for $z\in K^p,p=1,\ldots , m$,
$$
|f_{\z,l}(z)-t_{u,v,k,\z,l}(z)|\le 
\left\{\begin{array}{ll}
   \displaystyle
|f_{\z,l}(z)-p_{v,k,\z,K^l}(z)|+ |p_{v,k,\z,K^l}(z)|
|h_l(z)-q_{u,l}(z)|
 &\mb{ if } p=l,
\\[2ex]
|p_{v,k,\z,K^l}(z)|
|h_l(z)-q_{u,l}(z)|
&\mb{ if }p\neq l,
\end{array}\right.
$$
by (\ref{2.7}), (\ref{2.8}), and
the Bernstein-Walsh lemma 
(see \cite[p. 77]{wal} or \cite[p. 153]{saftot}),
we obtain the following estimates:

if $l=j$, then
\beq\label{2nn.11}
|f_{\z,l}(z)-t_{u,v,k,\z,l}(z)|\le 
\left\{\begin{array}{ll}
   \displaystyle
c_7\frac{d(\z,K^j)^k}{|\z-z|^{k+1}}+\frac{c_7 +1}{d(\z,K^j)}
e^{-ur_K}
 &\mb{ if } p=l,
\\[2ex]
 \displaystyle
\frac{c_7+1}{d(\z,K^j)}e^{vR_K-ur_K}
&\mb{ if }p\neq l;
\end{array}\right.
\eeq
\newpage
if $l\neq j$, then
$$
|f_{\z,l}(z)-t_{u,v,k,\z,l}(z)|
$$
\beq\label{2nn.12}
\le 
\left\{\begin{array}{lll}
   \displaystyle
	\frac{c_7}{n^k}+c_8e^{-ur_K}
 &\mb{ if } p=l,
\\[2ex]
\displaystyle
c_7\frac{d(\z,K^j)^k}{|\z-z|^{k+1}}+\frac{c_7 +1}{d(\z,K^j)}
e^{-ur_K}
 &\mb{ if } p=j,
\\[2ex]
 \displaystyle
c_8e^{vR_K-ur_K}
&\mb{ if }p\neq l, p\neq j.
\end{array}\right.
\eeq

Therefore, for the polynomial
$$
t_{u,v,k,\z}:=\sum_{l=1}^m t_{u,v,k,\z,l}\in\bP_{u+v}
$$
 according to (\ref{2nn.11}) and (\ref{2nn.12}) for $\z$ as in (\ref{2.7})
and $z\in K^p$, we obtain:

if $p=j$, then
\begin{eqnarray}
&&\left|\frac{1}{\z-z}-t_{u,v,k,\z}(z)\right|
\nonumber\\
\label{2nn.13}
&\le&
m\left(c_7\frac{d(\z,K^j)^k}{|\z-z|^{k+1}}+\frac{c_7 +1}{d(\z,K^j)}
e^{-ur_K}\right) ;
\end{eqnarray}

if $p\neq j$, then
\begin{eqnarray}
\left|\frac{1}{\z-z}-t_{u,v,k,\z}(z)\right|&\le&
\frac{c_7 +1}{d(\z,K^j)}
e^{vR_K-ur_K}\nonumber\\
\label{2nn.14}
&&+\frac{c_7}{n^k}+c_8e^{-ur_K}+(m-2)c_8e^{vR_K-ur_K}.
\end{eqnarray}
Let
$$
u:=\left\lfloor \frac{2R_K(n-1)}{2R_K+r_K}\right\rfloor,\quad
v:=\left\lfloor \frac{r_K(n-1)}{2R_K+r_K}\right\rfloor
$$
Note that $v\ge n/c_{9}$. To be sure that (\ref{2nn.13}) and (\ref{2nn.14}) hold
we need to have $c_3/v\le c^*/n$ which dictates the choice $c^*:=c_3c_{9}.$

Thus, using the L\"owner inequality (see \cite[p. 359, Corollary 2.5]{andbla}),
$
d(\z,K^j)\ge c_{10}/ n^2, c_{10}=c_{10}(K,k),
$
we obtain a polynomial
$$
s_{n-1,k,\z}:=t_{u,v,k,\z}\in\bP_{n-1}
$$
satisfying, 
by virtue of (\ref{2nn.13}) and (\ref{2nn.14}),
for $\z\in K^j_{r_K}\setminus K^j$ with 
$|\Phi_{\Om^j}(\z)|-1\ge c^*/n$, where
$c^*=c^*(K,k)$ and $n>n_1=n_1(K,k)$,
 the inequality
\beq\label{2.11}
\left|\frac{1}{\z-z}-s_{n-1,k,\z}(z)\right|\le c_{11}
\left\{\begin{array}{ll}
   \displaystyle
 \frac{d(\z,K^j)^k}{|\z-z|^{k+1}}
 &\mb{ if } z\in K^j,
\\[2ex]
n^{-k}&\mb{ if }z\in K\setminus K^j,
\end{array}\right.
\eeq
where $c_{11}=c_{11}(K,k)$.

\absatz{Chebyshev polynomials for a system of continua}

We start with the proof of the following estimate.

\begin{lem}\label{lem3.1}
Let $K$ be as in (\ref{1.1}). Then for $k\in\N$,  
\beq\label{3.1}
t_n(K)\le  c_1\sum_{j=1}^m\left|\left|\int_{L_{c^*/n,j}}\frac{d(\z,K^j)^{k}|d\z|}{|
\z-\cdot|^{k+1}}\right|\right|_{K^j},\quad n\ge n_1,
\eeq
where $c^*$ and $n_1$ are the constants from (\ref{2.11}), $c_1=c_1(K,k)$, and
$$
L_{\de,j}:=\{\z\in\Om^j:|\Phi_{\Om^j}(\z)|=1+\de\},\quad \de>0.
$$
\end{lem}
{\bf Proof.}
Let $F_n$ be defined by (\ref{2n.1}).
By our assumption
 $n$ is so large that the
curves $S_{n,j}:=L_{c^*/n,j}\subset\Om_{1/n^2}$ are mutually disjoint. Let
$S_n:=\cup_{j=1}^mS_{n,j}$. 

By \cite[p. 23, Lemma 2.3]{andbla}, for $\z\in\Om^j$, 
$w:=\Phi_{\Om^j}(\z)$, and $\Psi_{\Om^j}:=\Phi_{\Om^j}^{-1}$,
we have 
\beq\label{3n.3}
|\Psi_{\Om^j}'(w)|\asymp\frac{d(\z,K^j)}{|w|-1}.
\eeq
Therefore,
\begin{eqnarray}
|S_{n,j}|&=&\int_{|w|=1+c^*/n}|\Psi_{\Om^j}'(w)||dw|\nonumber\\
\label{3n.1}
&\asymp&\frac{n}{c^*}\int_{|w|=1+c^*/n}d(\Psi_{\Om^j}(w),K^j)|dw|
\le c n,\quad c=c(K,k).
\end{eqnarray}
By the Cauchy formula
$$
F_n(z)=\frac{1}{2\pi i}\sum_{j=1}^m\int_{S_{n,j}}\frac{W_n(\z)}{\z-z}d\z,\quad
z\in K.
$$
We can certainly assume that $k>1$.
Consider polynomial $F^*_n(z)=\al_nz^n+\ldots\in\bP_n$ defined as follows
$$
F_n^*(z):=\frac{1}{2\pi i}\sum_{j=1}^m\int_{S_{n,j}}W_n(\z)
\left(\frac{1}{\z-z}-s_{n-1,k,\z}(z)\right)d\z,\quad
z\in K,
$$
where $s_{n-1,k,\z}\in\bP_{n-1}$ satisfies (\ref{2.11}).

Since by (\ref{2.6}), for $\z\in S_n$,
$$
|W_n(\z)|\le|\Phi_\Om(\z)|^n=\exp(n g_\Om(\z))\le c_2,
$$
where $c_2=c_2(K,k)$, according to (\ref{2.11}) and (\ref{3n.1}), for $z\in K^j$, we obtain
\begin{eqnarray*}
|F^*_n(z)|&\le& c_3\left(\int_{S_{n,j}}\frac{d(\z,K^j)^k}{|\z-z|^{k+1}}|d\z|+
n^{-k}\sum_{l=1,l\neq j}^m |S_{n,l}|
\right)\\
&\le&
c_4\int_{S_{n,j}}\frac{d(\z,K^j)^k}{|\z-z|^{k+1}}|d\z|,
\end{eqnarray*}
where $c_l=c_l(K,k),l=3,4.$

Making use of
 (\ref{2.2}) and the obvious inequality
$t_n(K)\kap(K)^n\le||F_n^*||_K/|\al_n|$ we finally obtain (\ref{3.1}).

\hfill$\Box$

{\bf Proof of Theorem \ref{th1}.}
Changing the variable in the integrals from (\ref{3.1}) 
and using (\ref{3n.3}),
for sufficiently large $n$, we obtain
\begin{eqnarray}
&&\left|\left|\int_{L_{c^*/n,j}}\frac{d(\z,K^j)|d\z|}{|
\z-\cdot|^{2}}\right|\right|_{K^j}\nonumber\\
\label{3.2}
&\asymp& \frac{1}{n}
\left|\left|\int_{|w|=1+c^*/n}\left(\frac{|\Psi'_{\Om^j}(w)|}{|\Psi_{\Om^j}(w)-\cdot|}\right)^{2}
|dw|\right|\right|_{K^j}.
\end{eqnarray}
Furthermore, since by
 \cite[Chapter IX, \S 4, Lemma 3]{sue},
 $$
 \left|\left|\int_{|w|=1+c^*/n}\left(\frac{|\Psi'_{\Om^j}(w)|}{|\Psi_{\Om^j}(w)-\cdot|}\right)^{2}
|dw|\right|\right|_{K^j}\ole n\log n,
$$
the inequalities (\ref{3.1}) (with $k=1$) and (\ref{3.2}) imply  (\ref{1.2}).

\hfill$\Box$

Theorem \ref{th2} is a particular case of a more general result which we
 describe below. 
Let $K$ consist of one component, i.e., $m=1$, and
let $\Om$ be a John domain which can be defined as follows
 (see \cite[p. 98]{pom}).
 For a crosscut $\ga\subset\Om\setminus\{\infty\}$ of $\Om$
let  $H(\ga)$ be a bounded component of $\Om\setminus\ga$.
We say that $\ga$ is
  a  circular crosscut if
 $\ga\subset\Om\cap C(z,r)$ for some $z\in \partial\Om=\partial K, r>0$,
 and $z\in\ov{H(\ga)}$.
 Here  $C(z,r):=\{\z:|\z-z|=r\}$.
Then $\Om$ is a  John domain
if there exists a constant $\la_\Om>1$ such that
for any circular crosscut $\ga$ of $\Om$,
\beq\label{3.3}
\mb{diam}(H(\ga))\le \la_\Om |\ga|.
\eeq
By virtue of (\ref{1.3})
the complement of a quasiconformal arc as well as the unbounded
Jordan domain with a quasiconformal boundary both
  are John domains.

According to (\ref{3.3}) the function $\Psi_\Om$ has a continuous  extension to $\ov{\D^*}$
which we denote by the same letter $\Psi_\Om$. Next, we assume  that $\partial K$ is 
 piecewise
quasiconformal, i.e., there exist
$$
\theta_1<\theta_2<\ldots<\theta_{p}<\theta_{p+1}:=\theta_1+2\pi,
\quad p\ge2
$$
such that each $J_l:=\Psi_\Om(J'_l), l=1,\ldots,p$, where $J_l':=\{e^{i\theta}:
\theta_l\le\theta\le\theta_{l+1}\}$ is a quasiconformal arc.

Let
$$
z_l:=\Psi_\Om(e^{i\theta_j}),\quad \Ga'_l:=\{ re^{i\theta_j}:r\ge1\},\quad \Ga_l:=\Psi_\Om(\Ga_l'),
$$
$$
\Om_l':=\{re^{i\theta}:\theta_l<\theta<\theta_{l+1},r>1\},\quad \Om_l:=\Psi_\Om(\Om_l').
$$
By \cite[Lemma 2]{and15},
\beq\label{3.4}
|\z-z_l|\ole d(\z,K),\quad \z\in\Ga_l.
\eeq
Moreover, according to \cite[(4.14)]{and15},
\beq\label{3.5}
|\Ga_l(\z_1,\z_2)|\ole |\z_2-\z_1|,\quad \z_1,\z_2\in \Ga_l.
\eeq
Here for any arc or unbounded curve $\Ga\subset\C$ and $\z_1,\z_2\in\Ga$,
we denote by $\Ga(\z_1,\z_2)$ the bounded subarc of $\Ga$ between these points.

Thus, by virtue of (\ref{3.4}) and (\ref{3.5}), the curve $L_l^*:=\partial\Om_l
=\Ga_l\cup J_l\cup\Ga_{l+1}$ satisfies
$$
\mb{diam}(L_l^*(\z_1,\z_2))\ole |\z_2-\z_1|,\quad \z_1,\z_2\in L_l^*,
$$
i.e., by the Ahlfors criterion (see \cite[p. 100]{lehvir}), $L_l^*$ is quasiconformal.
Since by the same Ahlfors criterion $\partial\Om_l'=\Ga_l'\cup J_l'\cup\Ga'_{l+1}$
is also quasiconformal,
 the restriction of $\Phi_\Om$ to $\Om_l$ can be extended to a $Q_l$-quasiconformal
homeomorphism $\Phi_l:\C\to\C$ with some $Q_l\ge 1$ (see \cite[p. 98]{lehvir}).

The following result describes the distortion properties of
 $\Phi_l$ and  the inverse
mapping  $\Phi^{-1}_l$ which both are $Q$-quasiconformal
with $Q:=\max_{l=1,\ldots,p}Q_l$.
 \begin{lem}
\label{lem2.2} (\cite[p. 29]{andbla}) Let
 $F:\OC\to\OC$ be a $Q$-quasiconformal mapping, $Q\ge1$, with
  $F(\infty )=\infty$. Let $\z_k\in \C,w_k:=F(\z_k), k=1,2,3$,
  be such that $|w_1-w_2|\le c_5|w_1-w_3|$. Then $|\z_1-\z_2|\le
  c_6|\z_1-\z_3|$ and, in addition,
$$
\frac{1}{c_7}\left| \frac{w_1-w_3}{w_1-w_2} \right| ^{1/Q} \le
\left|\frac{\z_1-\z_3}{\z_1-\z_2}\right|\le
c_7\left|\frac{w_1-w_3}{w_1-w_2}\right|^{Q},
$$
where $c_j=c_j(c_5,Q),j=6,7$.
\end{lem}
We claim that for $z\in\partial K\setminus J_l$,
\beq\label{3.6}
d(z,J_l)\ole d(z,\Om_l).
\eeq
Indeed, let $z'_l\in\partial\Om_l$ be such that $|z-z'_l|=d(z,\Om_l)$.
The nontrivial  case arises when $z'_l\not\in J_l$, i.e., $z'_l\in\Ga_k$
for $k=l$ or $k=l+1$. Then by (\ref{3.4}) we obtain
$$
d(z,J_l)\le|z-z_k|\le |z-z'_l|+|z'_l-z_k|\ole |z-z'_l|=d(z,\Om_l),
$$
which yields (\ref{3.6}).

For $z\in \partial K$, denote by $z_l^*$ any point of $J_l$
with the property $|z-z_l^*|=d(z,J_l)$.
We claim that
\beq\label{3.7}
|\z-z_l^*|\ole |\z-z|,\quad \z\in\Om_l,\, z\in\partial K\setminus J_l.
\eeq
Indeed, by (\ref{3.6}),
$$
|\z-z_l^*|\le |\z-z|+|z-z_l^*|\ole|\z-z|+d(z,\Om_l)\le 2|\z-z|
$$
and (\ref{3.7}) follows.

For $\z\in\Om\setminus\{\infty\}$ denote by $\z_K:=\Psi_\Om(\Phi_\Om(\z)/|\Phi_\Om(\z)|)$
the ``projection" of $\z$ on $K$.
As an immediate application of Lemma \ref{lem2.2},
for $\z\in\Om_l$ and $z\in J_l$, we have
$$
\frac{d(\z,K)}{|\z-z|}\le\left|\frac{\z-\z_K}{\z-z}\right|
\ole \left(\frac{|\Phi_l(\z)|-1}{|\Phi_l(\z)-\Phi_l(z)|}\right)^{1/Q}.
$$
Now let $K$ be as in (\ref{1.1}).
We assume that each $\Om^j$ is a John domain
 and each  $\partial K^j$
is piecewise quasiconformal, i.e., each $\partial K^j$ consists of
$p_j$ quasiconformal arcs $J_{l,j},l=1,\ldots,p_j$
as described above.
 Let $\Phi_{l,j}$ be the appropriate
quasiconformal homeomorphism of $\C$ which is conformal in
a subdomain $\Om_l^j$ of $\Om^j$
with $J_{l,j}\subset\partial\Om^j_l$.
For $z\in \partial K^j$, denote by $z_{l,j}^* $ the nearest to $z$ point of
$J_{l,j}$ and let $w_{l,j}:=\Phi_{l,j}(z_{l,j}^*),
J_{l,j}':=\Phi_{l,j}(J_{l,j}).$

According to (\ref{3n.3}), Lemma \ref{lem2.2} with $F=\Psi_{\Om^j}$
restricted to $\Phi_{\Om^j}(\Om_l^j)$
and the triplet of points $\tau,\tau/|\tau|,w_{l,j}$, as well as (\ref{3.7}),
for $z\in\partial K^j$, $s=c^*/n$, and sufficiently large $n$, we obtain
\begin{eqnarray*}
\int_{L_{s,j}}\frac{d(\z,K)^k|d\z|}{|\z-z|^{k+1}}&=&
\sum_{l=1}^{p_j}\int_{L_{s,j}\cap\Om^j_l}\frac{d(\z,K)^k|d\z|}{|\z-z|^{k+1}}
\ole \sum_{l=1}^{p_j}
\int_{L_{s,j}\cap\Om^j_l}\frac{d(\z,K)^k|d\z|}{|\z-z_{l,j}^*|^{k+1}}\\
&\ole&
\sum_{l=1}^{p_j}
\frac{1}{s}\int_{|\tau|=1+s,\tau/|\tau|\in J'_{l,j}}\left|\frac{\Psi_{\Om^j}(\tau)-
\Psi_{\Om^j}(\tau/|\tau|)}{\Psi_{\Om^j}(\tau)-
\Psi_{\Om^j}(w_{l,j})}\right|^{k+1}|d\tau|\\
&\ole&
\sum_{l=1}^{p_j}
\frac{1}{s}\int_{|\tau|=1+s}\frac{s^{(k+1)/Q}|d\tau|}{|\tau-w_{l,j}|^{(k+1)/Q}}\ole
1
\end{eqnarray*}
if we fix  $k$ satisfying $k+1>Q$.
 
Comparing the last estimate with Lemma  \ref{lem3.1}  we obtain
the following statement.
\begin{th}\label{th4} 
Let $K$ be as in (\ref{1.1}). Assume that each $\Om^j$ is a John domain  
and each $\partial K^j$
 is piecewise
quasiconformal.
 Then (\ref{1.4}) holds.
\end{th}
This theorem yields Theorem \ref{th2}.

\absatz{Chebyshev polynomials for uniformly perfect sets}

We introduce some
definitions and notations from   geometric function theory. 
Let $K\subset\R$ be a
uniformly perfect set satisfying
\beq\label{4p.5}
 K\subset I:=[-1,1],\pm1\in K\neq I.
\eeq
  The open (with respect to
$\R$) set $I\setminus K$ consists of either a finite number $N\ge1$
or an infinite number $N=\infty$ of  disjoint open intervals, i.e.,
$$
I\setminus K=\bigcup_{j=1}^N(\al_j,\be_j),
$$
where
$(\al_j,\be_j)\cap(\al_k,\be_k)=\emptyset$ for $j\neq k$.

 It
follows immediately from (\ref{1.5}) that $\Om$ is regular
(for the Dirichlet problem), see
\cite{ran}, \cite{saftot}, i.e.,
  $g_\Om$
extends continuously to $K$ and $g_\Om(x):=0,\, x\in K$.
 Moreover, the Green function satisfies
\beq\label{4.1}
g_\Om(\z)\le c_1\, d(\z,K)^\al,\quad \z\in\C,
\eeq
where
 constants $c_1$ and $\al$ 
 could depend only on $\la_K$ from (\ref{1.5}),
see \cite[Lemma 4.1]{jerken} or \cite[p. 119]{garmar}.

We need the Levin conformal mapping which can be defined as follows
(for details, see
\cite{lev89}, \cite{and04}).
 Consider the  univalent in the upper half-plane $\He:=\{z:\Im z>0\}$ function
 $$
\phi(z)=\phi(z,K):=\pi+i\left(\int_K\log(z-\z)\,
d\mu(\z)-\log\kap(K)\right), \quad z\in \He ,$$
where $\mu=\mu_K$ is the equilibrium measure for  $K$.
 It  maps $\He$
 onto  a vertical half-strip
with $N$  slits parallel to the imaginary axis, i.e., the domain
 \beq\label{4.2}
 \Sigma_K:=\{ w:\, 0<\Re w<\pi,\Im w>0\}\setminus
\bigcup_{j=1}^{N}[u_j,u_j+iv_j],
\eeq
 where
$0<u_j=u_j(K)<\pi$ and $v_j=v_j(K)>0.$

 The continuous extension of $\phi$ to $\ov{\He}$ satisfies the following
boundary correspondence
$$\phi(\infty)=\infty,\, \phi((-\infty,-1])=\{w:\, \Re w=0,\Im w\ge
0\},$$
 $$ \phi([1,\infty))=\{w:\, \Re w=\pi,\Im w\ge 0\},\quad
  \phi(K)=[0,\pi],$$
 $$\phi([\al_j,\be_{j}])=[u_j,u_j+iv_j],\quad j=1,\ldots,N.$$
 Note
that in the last relation each point of
$[u_j,u_j+iv_j)$ has two preimages on
$[\al_j,\be_{j}]$.

The crucial fact is that $\phi$ satisfies
  \beq\label{4p.1}
g_\Om(z)=\Im\{\phi(z)\},\quad z\in \ov{\He}.
 \eeq
For a horizontal crosscut $\ga$ of $\Sigma_K$,
i.e., an interval
$\ga=(a+ib,c+ib)\subset\Sigma_k$ with endpoints on $\partial\Sigma_K$,
 denote by
$h(\ga)$ its ''height", that is,  $h(\ga):=b$.
\begin{lem}\label{lem4.1}
Any horizontal crosscut $\ga$ of $\Sigma_K$
with the property
$h(\ga)\le \sup_jv_j$ satisfies
\beq\label{4.3}
h(\ga)\le c_2|\ga|,\quad c_2=c_2(\la_K).
\eeq
\end{lem}
{\bf Proof.} For convenience, let $u_{-1}:=0$ and $u_0:=\pi$.
Let $\ga=(u_j+ih(\ga),u_k+ih(\ga))$ and $R:=\{w=u+iv:u_j<u<u_k,0<v<h(\ga)\}.$
 Denote by $\Ga'$
 the family of crosscuts of $\Sigma_K\cap R$
 which join $(u_j,u_k)$ to $\ga$ and let $\Ga^*$ be the family of crosscuts
 of the rectangle $ R$
 which join its horizontal boundary intervals. We refer to \cite{ahl}, \cite{lehvir}, \cite{garmar}
 for the basic properties of the module of a family of curves and arcs
 (such as conformal invariance, comparison principle, composition laws, etc.)
 We  use these properties without further citation.

 For the modules of $\Ga'$ and $\Ga^*$
we have
\beq\label{4.4}
m(\Ga')\le m(\Ga^*)=\frac{|\ga|}{h(\ga)}.
\eeq
At the same time, we claim that for the module of
$\Ga:=\phi^{-1}(\Ga')$ the estimate
\beq\label{4.5}
m(\Ga)\ge c_3,\quad c_3=c_3(\la_K)
\eeq
holds.

Indeed, without loss of generality, we assume that $j,k\ge1$ and
$\be_j-\al_j\le \be_k-\al_k$. 
The other particular cases may be handled in much the same way.
Denote by $\Ga_1$ the family of all crosscuts of
$$
G_1:=\{z=\al_j+re^{i\theta}:\be_j-\al_j<r<2(\be_j-\al_j),0<\theta<\pi\}
$$
which join $F_1:=K\cap[\be_j,2\be_j-\al_j]$ with $[3\al_j-2\be_j,2\al_j-\be_j]$.
Since $\Ga_1$ is ``fewer and longer" than $\Ga$, the comparison principle yields
\beq\label{4.6}
m(\Ga_1)\le m(\Ga).
\eeq
Note that by (\ref{1.5}),
\beq\label{4.7}
\kap(F_1)\ge\la_K(\be_j-\al_j).
\eeq
Consider the conformal mapping of $G_1$ onto
$$
G_2:=\{ w=re^{i\theta}:r_0<r<1,0<\theta<\pi\},\quad r_0:=\exp\left(-\frac{\pi^2}{\log 2}\right),
$$
given by the function
$$
w=f(z):=\exp\left(\frac{i\pi}{\log 2}\log\frac{z-\al_j}{\be_j-\al_j}\right)
$$
with the boundary correspondence
$$
f([\be_j,2\be_j-\al_j])=\{w=e^{i\theta}:0\le\theta\le\pi\},
$$
$$
f([3\al_j-2\be_j,2\al_j-\be_j])=\{w=r_0e^{i\theta}:0\le\theta\le\pi\}.
$$
Since for $\be_j\le x_1<x_2\le 2\be_j-\al_j$,
$$
|f(x_2)-f(x_1)|\ge \frac{x_2-x_1}{2(\be_j-\al_j)}\, ,
$$
by the Fekete-Szeg\H{o} Theorem (see \cite[p. 153]{ran})
and (\ref{4.7}) for the set $F_2:=f(F_1)$ we have
$$
\kap(F_2)\ge\frac{\kap(F_1)}{2(\be_j-\al_j)}\ge \frac{\la_K}{2}\, .
$$
Furthermore, let $\Ga_2:=f(\Ga_1)$ and denote by $\Ga_3$ the family
of all crosscuts of the annulus $\{\tau:r_0<|\tau|<1\}$ which join
$F_3:=F_2\cup \ov{F_2}$, where $\ov{F_2}:=\{\tau:\ov{\tau}\in F_2\}$, with
the circular boundary component $\{\tau:|\tau|=r_0\}$.
By the symmetry principle $m(\Ga_3)=2m(\Ga_2)$. Now we apply Pfluger's theorem
(see \cite[p. 212]{pom}) to obtain
\begin{eqnarray*}
m(\Ga_3)&\ge&\pi\left(\log\frac{1+r_0}{\sqrt{r_0}\kap(F_3)}\right)^{-1}\\
&\ge&\pi\left(\log\frac{1+r_0}{\sqrt{r_0}\kap(F_2)}\right)^{-1}\ge
\pi\left(\log\frac{2(1+r_0)}{\sqrt{r_0}\la_K}\right)^{-1}=:2c_3.
\end{eqnarray*}

Therefore, the conformal invariance of the module yields
$$
m(\Ga_1)=m(\Ga_2)=\frac{1}{2}m(\Ga_3)\ge c_3,
$$
which together with (\ref{4.6}) implies (\ref{4.5}).

At last, by virtue of the conformal invariance of the module, as well as
(\ref{4.4}) and (\ref{4.5}), we have
(\ref{4.3}) with $c_2:=c_3^{-1}$.

 \hfill$\Box$

Let now $1\le N<\infty$.
According to \cite{wid}, for $n\in\N$, either $\Phi_\Om(z)^n$ is single-valued
or it is multiple-valued. In the first case, we set $W_n(z):=\Phi_\Om(z)^n$
and in the second case there exist $q\le N$ points $x_{1,n},\ldots ,x_{q,n}\in I\setminus K$,
such that 
$$W_n(z):=\Phi_\Om(z)\prod_{l=1}^q\Phi_\Om(z,x_{l,n})^{-1},\quad z\in\Om,
$$
is single-valued in $\Om$.
According to \cite[pp. 159, 211]{wid}
each complementary interval 
$(\be_j, \al_{j+1})$ cannot have more than one point
from $\{x_{l,n}\}$.

Let polynomials $F_n(z)=F_n(z,K)$ be defined as in Section 2, i.e.,
$$
    F_n(z) := \frac{1}{2\pi i}  \int_{C_n}  \frac{W_n(\zeta)}
    {\zeta-z} \, d\zeta, \quad z\in\C,
  $$
 where $C_n\subset\Om\setminus\{\infty\}$ is a Jordan curve, 
oriented in the positive direction, containing
$K$ and $z$ in its interior.

By the Cauchy formula, for $z\in\Om\setminus\{\infty\}$ and sufficiently small
 $t>0$, we have
    $$
    F_n(z) = W_n(z) + \frac{1}{2\pi i} \; \int_{\tilde{K}_t} \; \frac{W_n(\zeta)}
    {\zeta-z} d\zeta,
    $$
    where $\tilde{K}_t:=\{\z\in\Om:d(\z,K)=t\}$ consists of $N+1$ disjoint curves
    each surrounding exactly one component of $K$.

Passing to the limit, we obtain for $z \in \Omega$ with $|z|<2$,
    \begin{eqnarray}
    |F_n(z)|&\le& |\Phi_\Om(z)|^n
    +\frac{1}{2\pi} \lim_{t \to 1^+} \; \int_{\tilde{K_t}}
    \frac{|\Phi_\Om(\zeta)|^n}{|\zeta-z|} \: |d\zeta| \nonumber\\
   & \le &
   e^{ng_\Om(z)}+\frac{1}{\pi}\int_{I\setminus D(z,d(z,K)}\frac{|d\z|}{|\z-z|}
   \nonumber\\
   \label{4.8}
&\ole& e^{ng_\Om(z)}+|\log d(z,K)|.
\end{eqnarray}
Here,
$$
D(z,r):=\{\z:|\z-z|<r\},\quad z\in\C,\, r>0.
$$
According to (\ref{4.1}) and (\ref{4.8}), for $z$ with the property $g_\Om(z)=1/n$,
we have the inequality
$$
|F_n(z)|\le c_4\log(n+1),\quad c_4=c_4(\la_K),
$$
which by the maximum principle for $F_n$ in $K_{1/n}$
is also  true for $z\in K$.

Note that $F_n(z)=\al_nz^n+\ldots$, where as in (\ref{2.2})
 \begin{eqnarray*}
 |\al_n|&=&
 \lim_{z\to\infty}\left|\frac{F_n(z)}{z^n}\right|=
\lim_{z\to\infty}\left|\frac{W_n(z)}{z^n}\right|\\
  &=&\kap(K)^{-n}\exp\left(-\sum_{l=1}^qg_\Om(x_{l,n})\right).
 \end{eqnarray*}
Therefore, by (\ref{4p.1}),
\beq\label{4.9}
t_n(K)\le\left|\left|\frac{F_n}{\al_n}\right|\right|_K\kap(K)^{-n}\le
c_4\log(n+1)\exp(V(K)),
\eeq
where
$$
V(K):=\sum_{j=1}^{N}v_j
$$
and $v_j=v_j(K)$ are defined by (\ref{4.2}).

{\bf Proof of Theorem \ref{th3}.}
Applying linear transformation if necessary we always can assume
 that
 $K$ satisfies (\ref{4p.5}).
By virtue of Theorem \ref{th2}, the only  nontrivial  case arises when $K$
consists of infinitely many components. 
Consider
$$
K^*_n:= I\cap \{z\in\C:g_\Om(z)\le1/n\},\quad n\in\N.
$$
It is worth pointing out that $K^*_n$ is uniformly perfect with $\la(K^*_n)=\la(K)$.
Moreover, by Lemma \ref{lem4.1}, $K_n^*$
consists of $N+1=N(K,n)+1\le c_5n$  disjoint closed intervals  and
\beq\label{4.10}
\kap(K)\le \kap(K_n^*)\le \kap(\{z\in\C:g_\Om(z)\le1/n\})=e^{1/n}\kap(K).
\eeq
Let $F_n(z)=F_n(z,K_n^*)$ be the Faber-Widom polynomial as above (constructed for $K_n^*$ instead
of $K$).
Denote by $v_{n,j}:=v_j(K_n^*),j=1,\ldots,N$, the quantities $v_j$ defined by
(\ref{4.2}) for $K_n^*$ instead of $K$. Note that
$$
\max_{1\le j\le N}v_{n,j}\le\sup_{1\le j<\infty}v_j(K)=c_{6}.
$$
For sufficiently large $n$, consider the sets
$$
\Lambda_0:=\left\{j:v_{n,j}\le\frac{1}{n}\right\},
$$
$$
\Lambda_k:=\left\{j:\frac{2^{k-1}}{n}< v_{n,j}\le\frac{2^{k}}{n}\right\},\quad
k=1,\ldots,k_0:=\lfloor \log_2(nc_{6})\rfloor +1.
$$
Since the number of elements in $\Lambda_0$ is at most $c_5n$ and by Lemma \ref{lem4.1}
the number of elements in
$
\Lambda_k$ is at most $ c_{7}n2^{-k},$ 
we obtain
$$
V(K_n^*)=\sum_{k=0}^{k_0}\sum_{j\in\Lambda_k}v_{n,j}\le c_5 +c_{7}k_0\le c_{8}\log n.
$$
Therefore,  by (\ref{4.9}) and (\ref{4.10}) 
$$
t_n(K)\le t_n(K^*_n)\frac{\kap(K^*_n)^n}{\kap(K)^n}
\ole n^{c_8}\log n,
$$
which implies (\ref{1.6}).

\hfill$\Box$

{\bf Acknowledgements.}
Part of this work was done during the Fall of 2016 semester, while the author visited
the Katholische Universit\"at Eichst\"att-Ingolstadt and
the Julius Maximilian University of W\"urzburg.
The author is  also grateful to   F. Nazarov 
 for his helpful comments.

\end{document}